\documentstyle[11pt,amssymb,amstex]{amsart}

\numberwithin{equation}{section}

\newtheorem{thm}{Theorem}[section]
\newtheorem{lem}[thm]{Lemma}
\newtheorem{cor}[thm]{Corollary}
\newtheorem{prop}[thm]{Proposition}
\newtheorem{defin}[thm]{Definition}
\newtheorem{rem}[thm]{Remark}
\newtheorem{prob}[thm]{Problem}
\newtheorem{obs}[thm]{Observation}
\newtheorem{exam}[thm]{Example}

\begin{document}

\title[Lotka-Volterra type operators]
{Stability and monotonicity of Lotka-Volterra type operators}

\author{Farrukh Mukhamedov}
\address{Farrukh Mukhamedov\\
Department of Computational \& Theoretical Sciences \\
Faculty of Sciences, International Islamic University Malaysia\\
P.O. Box, 141, 25710, Kuantan\\
Pahang, Malaysia} \email{{\tt far75m@@yandex.ru},{\tt
farrukh\_m@@iiu.edu.my}}

\author{Mansoor Saburov}
\address{Mansoor Saburov\\
Department of Computational \& Theoretical Sciences \\
Faculty of Science, International Islamic University Malaysia\\
P.O. Box, 141, 25710, Kuantan\\
Pahang, Malaysia} \email{{\tt msaburov@@gmail.com}}

\begin{abstract}
In the present paper, we study Lotka-Volterra (LV) type operators
defined in finite dimensional simplex. We prove that any LV type
operator is a surjection of the simplex. After, we introduce a new
class of LV-type operators, called $M$LV type. We prove convergence
of their trajectories and study certain its properties. Moreover, we
show that such kind of operators have totaly different behavior than
${\mathbf{f}}$-monotone LV type operators. \vskip 0.3cm \noindent
{\it
Mathematics Subject Classification}: 15A51, 47H60, 46T05, 92B99.\\
{\it Key words}: Lotka-Volterra type operators; stability; monotone
operator; simplex.
\end{abstract}

\maketitle

\section{Introduction}

Lotka-Volterra (LV) systems typically model the time evolution of
conflicting species in biology \cite{B,V1}. They have been largely
studied starting with Lotka \cite{L} and Volterra \cite{V2}. There
are many other natural phenomena modelled by LV systems (see
\cite{Ta}). On the other hand, the use of LV discrete-time systems
is a well-known subject of applied mathematics \cite{L1}. They were
first introduced in a biomathematical context by Moran \cite{Moran},
and later popularized by May and collaborators \cite{May,MO}. Since
then, LV systems have proved to be a rich source of analysis for the
investigation of dynamical properties and modelling in different
domains, not only population dynamics \cite{FG,BF,HJ,LW,NSE},  but
also physics \cite{PL,UR}, economy \cite{D}, mathematics
\cite{G,HS,L1,Ta,V,U}. Typically in all these applications, the LV
systems are taken quadratic. It is natural to investigate
non-quadratic LV systems. In \cite{GMM} were introduced,
generalizing the LV systems, to model the interaction among
biochemical populations. Cubic polynomial LV vector fields have
appeared explicitly modelling certain phenomena arising in
oscillating chemical reactions as the socalled
Lotka-Volterra-Brusselator (see \cite{FNSS}) and in well-known
predator-prey models that give rise to periodic variations in the
populations (see \cite{DV,LK}, etc). In \cite{GR} the global phase
portraits in the Poincare disc of the cubic polynomial vector fields
of LV type having a rational first integral of degree 2 is
classified. There, the linearizability problem for the
two-dimensional planar the cubic polynomial vector fields of LV type
having a rational first integral of degree 2, is investigated.
 The necessary and sufficient conditions for the linearizability of this system are found.

Recently, the family of discrete-time systems termed quasipolynomial
(QP) has attracted some attention in the literature \cite{HB1}. In
this context, it is worth noting that the interest of QP discrete
systems arises from several different features. In the first place,
they constitute a wide generalization of LV models. However, LV
systems are not just a particular case of QP ones but play a
central, in fact canonical role in the QP framework, as will be
appreciated in what is to follow. In \cite{HB2} the quasipolynomial
(QP) generalization of LV discrete-time systems is considered. Use
of the QP formalism is made for the investigation of various global
dynamical properties of QP discrete-time systems including
permanence, attractivity, dissipativity and chaos. The results
obtained generalize previously known criteria for discrete LV
models.

In \cite{Mu2} it is established new sufficient  conditions for
global asymptotic stability of the positive equilibrium in some
LV-type discrete models. Applying the former results
\cite{Mu1} on sufficient conditions for the persistence of
nonautonomous discrete LV systems, conditions for the persistence of
the above autonomous system is obtained, and extending a similar
technique to use a nonnegative Lyapunov like function offered by
\cite{SHM}, new conditions for global asymptotic stability of the
positive equilibrium is found.

The mentioned papers show importance the study of limiting behavior
of discrete LV type operators. Therefore, in \cite{GS}
${\mathbf{f}}$-monotone LV type operators on the simplex have been
defined. It was proved the existence of Lyapunov functions for such
operators which allowed to study limiting behaviors ones. Continuing
the previous investigations, in the present paper, we first show
that any LV type operator is a surjection of the simplex. After, we
introduce a new class of LV type operators, called $M$LV type, on
the simplex. We shall prove convergence of their trajectories and
study certain its properties. Moreover, we show that such kind of
operators have totaly different behavior than
${\mathbf{f}}$-monotone LV type operators.

\section{Preliminaries}

Let $$S^{m-1}=\left\{x=(x_1,x_2,\dots,x_m)\in{\mathbb{R}}^m:
\sum\limits_{k=1}^mx_k=1, x_k\geq0\right\}$$ be the
$(m-1)-$dimensional simplex. One can see that the points
$e_k=(\delta_{1k},\delta_{2k},\dots,\delta_{mk})$ are  the extremal
points of the simplex $S^{m-1},$ where $\delta_{ik}$ is the
Kronecker's symbol.

Let $I=\{1,2,\dots,m \}$ and $\alpha $  be an arbitrary subset of
$I$. The set
$$\Gamma_\alpha= \{x\in S^{m-1}: x_k=0, \ k\notin \alpha
\}$$ is called \emph{a face} of the simplex. \emph{A relatively
interior} $ri\Gamma_\alpha$ of the face $\Gamma_\alpha$ is defined
by
$$ri\Gamma_\alpha= \{x\in \Gamma_\alpha: x_k>0,
k\in\alpha \}.$$ \emph{The center} of the face $\Gamma_\alpha$ is
defined by
$$\bigg(0,\dots,0,\underbrace{\frac{1}{|\alpha|}}_{i_1},0,\dots,0,\underbrace{\frac{1}{|\alpha|}}_{i_2},0,\dots,0,\underbrace{\frac{1}{|\alpha|}}_{i_r},0,\dots,0\bigg)$$
where $\alpha=\{i_1,i_2,\dots,i_r\},$ $i_1<i_2<\cdots<i_r$ and
$|\alpha|$ is the cardinality of a set $\alpha.$

Given a mapping ${\textbf{f}}:x\in S^{m-1}\to
(f_1(x),f_2(x),\dots,f_m(x))\in{\mathbb{R}}^m$ in what follows, we
are interested in the following operator defined by
\begin{eqnarray}\label{Volterra}
(Vx)_k=x_k(1+f_k(x)), \ k=\overline{1,m} \ \ x\in S^{m-1}.
\end{eqnarray}

\begin{prop}\label{defin}
Let  $V$ be an operator given by \eqref{Volterra}. The following
conditions are equivalent:
\begin{itemize}
  \item [(i)] The operator $V$ is continuous in $S^{m-1}$ and
  $V(S^{m-1})\subset S^{m-1}.$
  Moreover, $V(ri\Gamma_\alpha)\subset ri\Gamma_\alpha$ for all $\alpha\subset I.$
  \item [(ii)]The mapping  ${\emph{\textbf{f}}}\equiv(f_1,f_2,\dots,f_m):S^{m-1}\rightarrow
{\mathbb{R}}^m$ satisfies the following conditions:
\begin{itemize}
  \item [$1^0$] $\emph{{\textbf{f}}}$ is continuous in $S^{m-1};$
  \item [$2^0$] for every $x\in S^{m-1}$ one has $f_k(x)\ge-1,$   for all $k=\overline{1,m};$
  \item [$3^0$] for every $x\in S^{m-1}$ one has $\sum\limits_{k=1}^mx_kf_k(x)=0;$
  \item [$4^0$] for every $\alpha\subset I$ one holds $f_k(x)>-1$ for
  all $x\in ri\Gamma_\alpha$ and $k\in \alpha.$
\end{itemize}
\end{itemize}
\end{prop}

\begin{pf}
$(i)\Rightarrow(ii).$ The continuity of $V$ implies $1^0.$ Take
$x\in S^{m-1}$ and it yields that (a) $(Vx)_k\ge 0$; (b) $
\sum\limits_{k=1}^m(Vx)_k=1$. Hence, from (a) it follows that
$x_k(1+f_k(x))\ge0$ which implies $2^0.$ From (b) one has
$$\sum\limits_{k=1}^mx_k+\sum\limits_{k=1}^mx_kf_k(x)=1$$ which
immediately yields $3^0.$

Let $x\in ri\Gamma_\alpha$, then $Vx\in ri\Gamma_\alpha$ which with
\eqref{Volterra} and $x_k>0$ for all $k\in\alpha$ implies that
$f_k(x)>-1$ for all $k\in\alpha$ this means $4^0.$

The implication $(ii)\Rightarrow(i)$ is evident.
\end{pf}

We say that an operator $V$ defined by \eqref{Volterra} is
\emph{Lotka-Volterra (LV) type} if one of the conditions of
Proposition \ref{defin} is satisfied. The corresponding mapping
$\textbf{f}$ is called \emph{generating mapping} for $V.$  From
Proposition \ref{defin} we immediately infer that any LV type
operator maps the simplex $S^{m-1}$ into itself. By ${\mathcal{V}}$
we denote the set of all LV type operators.

Note that first such kind of operators were considered in \cite{GS}
and there were proved the following theorem.

\begin{thm}[\cite{GS}]\label{linear} Let $\textbf{f}\equiv(f_1,f_2,\dots,f_m):{{\mathbb{R}}^m}\rightarrow
{{\mathbb{R}}^m}$ be a linear mapping. Then $\textbf{{f}}$ satisfies
the conditions $1^0-4^0$ if and only if one has
\begin{eqnarray}\label{LO}
f_k(x)=\sum\limits_{i=1}^ma_{ki}x_i, \ \ \ k=\overline{1,m}
\end{eqnarray}
with
$$a_{ki}=-a_{ik}, \ \ |a_{ki}|\le1 \ \ \ \forall \
k,i=\overline{1,m}.$$
\end{thm}

\begin{rem}\label{quadvolterra} Note that LV type operator $V:S^{m-1}\rightarrow
S^{m-1}$ with generating mapping given by \eqref{LO}, i.e.
\begin{eqnarray}\label{QSO}
(Vx)_k=x_k\bigg(1+\sum\limits_{i=1}^ma_{ki}x_i\bigg), \ \ \
k=\overline{1,m}
\end{eqnarray}
where $a_{ki}=-a_{ik}$, $|a_{ki}|\le1,$ for all
$k,i=\overline{1,m},$ is called quadratic volterrian operators. By
${\mathcal{QV}}$ we denote the set of all quadratic volterrian
operators. In \cite{G,MS,V,Z} limiting behavior of quadratic
volterrian operators were studied.
\end{rem}


Given $x^0\in S^{m-1},$ then the sequence
$\left\{x^0,Vx^0,V^{2}x^0,\cdots,V^{n}x^0,\cdots\right\}$ is called
\emph{a trajectory of $V$ starting from the point $x^0,$} where
$V^{n+1}x^0=V(V^{n}x^0), n=1,2,\dots$ By $\omega(x^0)$ we denote the
set of all limiting points of such a trajectory.

A point $x\in S^{m-1}$ is called \emph{fixed} if $Vx=x$ and by
$Fix(V)$ we denote the set of all fixed points of $V.$ A point $x\in
S^{m-1}$ is called \emph{$r$-periodic} if $V^rx=x$ and $V^ix\neq x$
for all $i\in\overline{1,r-1}.$

Let us introduce some necessary notations taken from \cite{SE} (see
also \cite{N,T}). Let $U$ be a bounded open subset of
${\mathbb{R}}^m.$ For $\textbf{f}\in C^1(U)$ the Jacobi matrix of
$\textbf{f}$ at $x\in U$ is
$\textbf{f}{'}(x)=(\partial_{x_i}f_j(x))_{i,j=1}^m$ and the Jacobi
determinant of $\textbf{f}$ at $x\in U$ is
$$J_{\textbf{f}}(x)=\det \textbf{f}{'}(x).$$ Given $y\in
{\mathbb{R}}^m$ let us set
\begin{eqnarray*}
D^r_y(\overline{U})=\{\textbf{f}\in
C^r(\overline{U}): y\notin \textbf{f}(\partial U)\}, \ \ \ D_y(\overline{U})=\{\textbf{f}\in
C(\overline{U}): y\notin \textbf{f}(\partial U)\}.
\end{eqnarray*}

A function $\textbf{deg}:\textbf{f}\in D_y(\overline{U})\to
\textbf{deg}(\textbf{f},U,y)\in {\mathbb{R}}$ is called \emph{degree
of \emph{\textbf{f}} at $y$}  if it satisfies the following
conditions:
\begin{itemize}
  \item [(D1)] $\textbf{deg}(\textbf{f},U,y)=\textbf{deg}(\textbf{f}-y,U,0)$;
  \item [(D2)] $\textbf{deg}(Id,U,y)=1$ if $y\in U$, where $Id$ is an identity mapping.
  \item [(D3)] If $U_{1},U_2$ are open, disjoint subsets of $U$ such
  that $y\notin \textbf{f}(\overline{U}\setminus(U_1\cup U_2)),$
  then
  $\textbf{deg}(\textbf{f},U,y)=\textbf{deg}(\textbf{f},U_1,y)+\textbf{deg}(\textbf{f},U_2,y)$;

  \item [(D4)] If $H(t)=(1-t)\textbf{f}+t\textbf{g}\in
  D_y(\overline{U}),$ $t\in[0,1],$ then
  $\textbf{deg}(\textbf{f},U,y)=\textbf{deg}(\textbf{g},U,y)$.
\end{itemize}

\begin{thm}[\cite{T}]\label{degree} There is a unique degree
\emph{\textbf{deg}}  satisfying (D1)-(D4). Moreover,
$\emph{\textbf{deg}}(\cdot,U,y):D_y(\overline{U})\to {\mathbb{Z}}$
is constant on each component and given $\emph{\textbf{f}}\in
D_y(\overline{U})$ we have
$$\emph{\textbf{deg}}(\emph{\textbf{f}},U,y)=\sum\limits_{x\in {\widetilde{\emph{\textbf{f}}}}(y)} sign
J_{\widetilde{\emph{\textbf{f}}}}(x)$$ where
$\widetilde{\emph{\textbf{f}}}\in D_y^2(\overline{U})$ is in the
same component of $D_y(\overline{U}),$ say
$\|\emph{\textbf{f}}-\widetilde{\emph{\textbf{f}}}\|<
\textbf{dist}(y,\emph{\textbf{f}}(\partial U)),$ such that $y\in
RV(\emph{\textbf{f}})=\{y\in {\mathbb{R}}^m| \forall x\in
\emph{\textbf{f}}^{-1}(y): J_{\emph{\textbf{f}}}(x)\neq0\}.$
\end{thm}

\begin{thm}[\cite{K}]\label{nonnulldeg}
If the degree of the mapping  $\textbf{\emph{f}}:S^{m-1}\to S^{m-1}$
is not zero then the mapping $\textbf{\emph{f}}$ is \textbf{onto}
(i.e. surjective).
\end{thm}

\section{Surjectivity  of Lotka-Volterra type operators}

In this section we shall show that all LV type operators are
surjective.

\begin{thm}\label{surVolterra} Any
LV type operator given by \eqref{Volterra} maps simplex $S^{m-1}$
onto itself. Namely, $V$ is a surjection of $S^{m-1}$.
\end{thm}
\begin{pf}
Consider a family od operators $V_\varepsilon:S^{m-1}\rightarrow
S^{m-1}$ given by
$$(V_\varepsilon x)_k=x_k(1+\varepsilon f_k(x)), \ \ k=\overline{1,m},$$
where $0\le\varepsilon\le 1,$ which homotopical connects an identity
mapping   $Id:S^{m-1}\rightarrow S^{m-1}$ and  the LV type operator
\eqref{Volterra}, i.e. $V_\varepsilon=(1-\varepsilon)Id +\varepsilon
V.$
 According to Theorem \ref{degree}, (D2) and (D4) we have $\mathrm{deg}(V)
 =\mathrm{deg}(Id)=1.$
Therefore, Theorem \ref{nonnulldeg} implies that LV type operator
\eqref{Volterra} is a surjection of $S^{m-1}$.
\end{pf}

It is well known (see \cite{K}) that any continuous bijective
mapping of a compact set to itself is homeomorphism of compact,
hence we have the following

\begin{cor}\label{gomeom}
Any LV type operator given by \eqref{Volterra} is homeomorphism of
the simplex if and only if it is  injective.
\end{cor}

\begin{rem}\label{exquadvolterra} Note that any quadratic
volterrian operator given by \eqref{QSO} is a homeomorphism of the
simplex (see \cite{G}). It is worth to note that not all LV type
operators are homeomorphisms of the simplex (see Example
\ref{notinjecperiodic}).
\end{rem}

\section{$M-$Lotka-Volterra type operators.}

In this section we introduce a class of LV type operators, called
$M-$LV, and study their asymptotic behavior.

Given $x\in S^{m-1}$ put $$M(x)=\{i\in I:
x_i=\max\limits_{k=\overline{1,m}}x_k\},$$ here as before
$I=\{1,\dots,m\}$.

\begin{defin} An LV type operator given by \eqref{Volterra} is
called $M_1-$Lotka-Volterra (for shortness $M_1$LV) (resp.
$M_0-$Lotka-Volterra ($M_0$LV)) if for each $x\in S^{m-1}$ and for
all $k\in M(x),$ $j=\overline{1,m}$ the functional
$$\varphi(x)=x_k-x_j$$ is increasing (res. decreasing)
along the trajectory of $V$ starting from the point $x,$ i.e.
$\varphi(V^{k}x)\le\varphi(V^{k+1}x), \ \ k\geq 0$ (resp.
$\varphi(V^{k}x)\ge\varphi(V^{k+1}x), \ \ k\geq 0)$\end{defin}

By ${{\mathcal{VM}}_1}$ and ${{\mathcal{VM}}_0}$ we denote the sets
of all $M_1$LV and $M_0$LV type operators, respectively.

Note that in \cite{GS} LV type operators with functionals of the
form $\varphi(x)=\prod_{k=1}^mx_k^{p_k}$ has been investigated.

\begin{rem}\label{interMvolterra} It immediately follows from the definition  that $${{\mathcal{VM}}_1}\cap{{\mathcal{VM}}_0}=\{Id\},$$
where $Id:S^{m-1}\to S^{m-1}$ is an identity mapping.
\end{rem}

\begin{prop}\label{convexkompo} Let $V_0$ and $V_1$ be $M_1$LV
(resp. $M_0$LV) type operators. Then the following conditions are
satisfied:
\begin{itemize}
  \item [(i)] The operator $V_1\circ V_0$ is $M_1$LV
(res. $M_0$LV) type.
  \item [(ii)] For each $\lambda\in[0,1]$ the operator $(1-\lambda) V_0+\lambda
  V_1$ is $M_1$LV (res. $M_0$LV) type.
\end{itemize}
\end{prop}

\begin{pf}
Without loss of generality we may suppose that the operator  $V_0$
and $V_1$ are $M_1$LV type. Then  for each $x\in S^{m-1}$ and for
all $k\in M(x),$ $j=\overline{1,m}$ we have
$$x_k-x_j\le(V_0x)_k-(V_0x)_j\le(V_1(V_0x))_k-(V_1(V_0x))_j$$
which implies that $V_1\circ V_2\in {\mathcal{VM}}_1.$

Now for all $\lambda\in[0,1]$ one finds
\begin{eqnarray*}
 x_k-x_j&=&(1-\lambda)(x_k-x_j)+\lambda(x_k-x_j)\\
 &\le& (1-\lambda)((V_0x)_k-(V_0x)_j)+\lambda((V_1x)_k-(V_1x)_j)\\
 &=&((1-\lambda)V_0x
 +\lambda V_1x)_k-((1-\lambda)V_0x
 +\lambda V_1x)_j,
\end{eqnarray*}
that yields the required assertion.

By the similar argument one can prove the statements for the case of
$M_0$LV type operators.
\end{pf}

\begin{cor}
The sets ${\mathcal{VM}}_1$ and ${\mathcal{VM}}_0$ are convex.
\end{cor}

Let us provide some  examples of $M_1$LV and $M_0$LV type operators,
respectively.

\begin{exam}\label{M1volterra} Let us consider an operator
$V_{\varepsilon,\ell}$ defined by
\begin{eqnarray}\label{VLLEPSILON}
(V_{\varepsilon,\ell}
x)_k=x_k\bigg(1+\varepsilon\bigg(x_k^\ell-\sum\limits_{i=1}^mx_i^{\ell+1}\bigg)\bigg),
\ k=\overline{1,m}
\end{eqnarray} where $0<\varepsilon\le 1$ and
$\ell\in {\mathbb{N}}.$
\end{exam}

Let us first show that operator given by \eqref{VLLEPSILON} is an LV
type.  One can immediately see that the generating mapping
$\textbf{f}_{\varepsilon,\ell}:S^{m-1}\to {\mathbb{R}}^m$ of
$V_{\varepsilon,\ell}$ is given by
$$\textbf{f}_{\varepsilon,\ell}(x)=\left(\varepsilon\bigg(x_1^\ell-\sum\limits_{i=1}^mx_i^{\ell+1}\bigg), \varepsilon\bigg(x_2^\ell-\sum\limits_{i=1}^mx_i^{\ell+1}\bigg),\cdots,
\varepsilon\bigg(x_m^\ell-\sum\limits_{i=1}^mx_i^{\ell+1}\bigg)\right)$$
which is obviously continuous. The inequality
\begin{eqnarray}\label{FLEPSILON}
(\textbf{f}_{\varepsilon,\ell}(x))_k+1&=&
\varepsilon\bigg(x_k^\ell+\sum\limits_{i=1}^mx_i(1-x_i^\ell)\bigg)+1-\varepsilon\ge0
\end{eqnarray}
shows that $(\textbf{f}_{\varepsilon,\ell}(x))_k\ge -1,$ for any
$x\in S^{m-1}$ and $k=\overline{1,m}.$

One can see that
\begin{eqnarray*}
\sum\limits_{k=1}^mx_k(\textbf{f}_{\varepsilon,\ell}(x))_k=\varepsilon
\bigg(\sum\limits_{k=1}^mx_k^{\ell+1}-\sum\limits_{i=1}^mx_i^{\ell+1}\sum\limits_{k=1}^mx_k\bigg)=0.
\end{eqnarray*}

Take any subset $\alpha\subset I$, then the inequality
\eqref{FLEPSILON} implies that
$(\textbf{f}_{\varepsilon,\ell}(x))_k>-1,$ for every $x\in
\mathrm{ri}\Gamma_\alpha$ and $k\in\alpha$. This means that
$V_{\varepsilon,\ell}$ is a Lotka-Volterra type operator.

Now let us establish  that  $V_{\varepsilon,\ell}$ is $M_1$LV type.
Indeed, take any $x\in S^{m-1}$ and one easily finds that
\begin{eqnarray}\label{E143}
(V_{\varepsilon,\ell} x)_k-(V_{\varepsilon,\ell}
x)_j=(x_k-x_j)\bigg(1+\varepsilon\sum\limits_{r=0}^\ell
x_k^{\ell-r}x_j^r-\varepsilon\sum\limits_{i=1}^mx_i^{\ell+1}\bigg).
\end{eqnarray}
On the other hand, the following relation
\begin{eqnarray}\label{E144}
1+\varepsilon\sum\limits_{r=0}^lx_k^{\ell-r}x_j^r-\varepsilon\sum\limits_{i=1}^mx_i^{\ell+1}=
1-\varepsilon+\varepsilon\sum\limits_{r=0}^\ell
x_k^{\ell-r}x_j^r+\varepsilon\sum\limits_{i=1}^mx_i(1-x_i^\ell)\ge0
\end{eqnarray}
with \eqref{E143} implies that
$$sign((V_{\varepsilon,\ell} x)_k-(V_{\varepsilon,l}
x)_j)=sign(x_k-x_j),$$ which means
$$M(V_{\varepsilon,\ell}x)=M(x).$$
Take $k\in M(x),$ then we have
\begin{eqnarray}\label{E145}
\varepsilon\sum\limits_{r=0}^\ell
x_k^{\ell-r}x_j^r-\varepsilon\sum\limits_{i=1}^mx_i^{\ell+1}=\varepsilon\sum\limits_{i=1}^mx_i(x_k^\ell-x_i^\ell)+\varepsilon\sum\limits_{r=1}^\ell
x_k^{\ell-r}x_j^r\ge0,
\end{eqnarray} for all $j=\overline{1,m}.$ Hence, from \eqref{E143}, \eqref{E145} we obtain
$$x_k-x_j\le(V_{\varepsilon,\ell} x)_k-(V_{\varepsilon,\ell}
x)_j$$ for every $k\in M(x)$ and $j\in I$. \\

By the similar argument used in Example \ref{M1volterra}, we can
show that the following operator $W_{\varepsilon,\ell}$ is $M_0$LV
type.

\begin{exam}\label{M0Volterra} Let us consider an operator $W_{\varepsilon,\ell}: S^{m-1}\rightarrow
S^{m-1}$ defined by
\begin{eqnarray}\label{WEPSILONL}
(W_{\varepsilon,\ell}
x)_k=x_k\left(1+\varepsilon\left(\sum\limits_{i=1}^mx_i^{\ell+1}-x_k^\ell\right)\right),
k=\overline{1,m}
\end{eqnarray} where $0<\varepsilon\le 1$ and $\ell\in {\mathbb{N}}.$
\end{exam}

Observe that by means of the provided examples and Proposition
\ref{convexkompo} one can construct lots of nontrivial examples of
$M_1$LV and $M_0$LV type operators, respectively.

 To study stability properties of $M_0$LV and
$M_1$LV type operators we need the following auxiliary result.

\begin{lem}\label{converge} If for a sequence $\{x^{(n)}\}_{n=0}^\infty\subset S^{m-1}$
and some $k\in I$ the limits
\begin{eqnarray}\label{Lem47}
\lim\limits_{n\to\infty}\left(x_k^{(n)}-x_j^{(n)}\right), \ \
\forall  \ \  j=\overline{1,m},
\end{eqnarray} exist, where
$x^{(n)}=(x_1^{(n)},x_2^{(n)},\dots,x_m^{(n)}),$ then the sequence
$\{x^{(n)}\}_{n=0}^\infty$ converges.
\end{lem}

\begin{pf}
The convergence of the sequences
$\left\{x_k^{(n)}-x_j^{(n)}\right\}_{n=0}^\infty,$ for all
$j=\overline{1,m},$ implies  the convergence of a sequence
$\bigg\{\sum\limits_{j=1}^m\big(x_k^{(n)}-x_j^{(n)}\big)\bigg\}_{n=0}^\infty.$
Then the equality
$$mx_k^{(n)}=\sum\limits_{j=1}^m\left(x_k^{(n)}-x_j^{(n)}\right)+\sum\limits_{j=1}^mx_j^{(n)}=\sum\limits_{j=1}^m\left(x_k^{(n)}-x_j^{(n)}\right)+1.$$
implies the convergence of the sequence
$\{x_k^{(n)}\}_{n=0}^\infty.$ From \eqref{Lem47} we obtain the
convergence of $\{x_j^{(n)}\}_{n=0}^\infty,$ for all
$j=\overline{1,m}$ which yields the convergence of
$\{x^{(n)}\}_{n=0}^\infty.$
\end{pf}

Now we are ready to prove stability property of $M_0$LV and $M_1$LV
type operators.

\begin{thm}\label{Volconverge}
Let $V$ be a $M_1$LV (resp. $M_0$LV) type operator. Then the
trajectory $\{V^nx\}_{n=0}^\infty$ converges for every $x\in
S^{m-1},$ i.e. $\omega(x)$ is a single point and $\omega(x)\in
Fix(V).$
\end{thm}

\begin{pf}
Let $V$ be a $M_1$LV type operator. Then for some $k\in M(x)$ and
all $j=\overline{1,m}$ we have
$$x_k-x_j\le(Vx)_k-(Vx)_j\le\cdot\cdot\cdot\le(V^nx)_k-(V^nx)_j\le\cdots\le1$$
Therefore the sequence
$\left\{(V^nx)_k-(V^nx)_j\right\}_{n=0}^\infty$ converges. It
follows from Lemma \ref{converge} that the trajectory
$\{V^nx\}_{n=0}^\infty$ converges.

By the similar way the statement can be proved for a $M_0$LV case.
\end{pf}

\begin{lem}\label{savemax}
Let $V$ be a $M_1$LV type operator. Then for every $x\in S^{m-1}$
and for all $n\in {\mathbb{N}}$ one has
$$M(V^nx)=M(x).$$ Moreover, if there exists a limit
$\lim\limits_{n\to\infty}V^nx=x^{*}$ then $M(x^{*})=M(x).$
\end{lem}

\begin{pf}
Let $k\in M(x).$ Since $V$ is $M_1$LV type then for every
$j=\overline{1,m}$ one has
\begin{eqnarray}\label{Lem48} 0\le x_k-x_j\le (Vx)_k -
(Vx)_j\le\cdot\cdot\cdot\le(V^nx)_k-(V^nx)_j\le\cdot\cdot\cdot
\end{eqnarray}
which implies that $k\in M(V^nx)$ i.e.
\begin{eqnarray}\label{Lem49} M(x)\subset M(V^nx).
\end{eqnarray} Now we show the inclusion
$M(V^nx)\subset M(x).$ Assume from the contrary, i.e.  there is
$k_0\in M(V^nx)$ that $k_0\notin M(x).$ Now take any  $k_1\in M(x)$,
then from \eqref{Lem49} we infer that $k_1\in M(V^nx)$ which means
$(V^nx)_{k_1}-(V^nx)_{k_0}=0$. On the other hand, from
\eqref{Lem48},\eqref{Lem49} we have that
$$0<x_{k_1}-x_{k_0}\le(Vx)_{k_1} - (Vx)_{k_0}\le\cdots\le(V^nx)_{k_1}-(V^nx)_{k_0}=0$$
which is a  contradiction, hence $M(V^nx)\subset M(x).$ Thus, we
have $M(V^nx)=M(x)$, for any $n\in {\mathbb{N}}.$

Now assume that $\{V^nx\}_{n=0}^\infty$ converges to $x^{*}.$ Then
from \eqref{Lem48} one has
\begin{eqnarray}\label{Lem410} x_k-x_j\le x_k^{*}-x_j^{*}
\end{eqnarray} for all $k\in
M(x)$ and $j=\overline{1,m}.$ Then \eqref{Lem410} yields that
\begin{eqnarray}\label{Lem411} M(x)\subset M(x^{*}).\end{eqnarray} Now we show the inverse
inclusion $M(x^{*})\subset M(x).$ Assume from the contrary, i.e.
there is $k_0\in M(x^{*})$ that $k_0\notin M(x).$ Then we use the
same argument as above, i.e. for any $k_1\in M(x)$ it follows from
\eqref{Lem410}, \eqref{Lem411} that
$$0<x_{k_1}-x_{k_0}\le x^{*}_{k_1}-x^{*}_{k_0}=0.$$ Again,
the last contradiction shows that $M(x^{*})\subset M(x)$ or  $M(x^{*})=
M(x).$
\end{pf}

\begin{rem}
Note that in general a similar result as Lemma \ref{savemax} is not
satisfied for  $M_0$LV type operators (see Observation
\ref{omegaWEPL}).
\end{rem}

\begin{thm}\label{FixedVol}
Let $V$ be an $M_1$LV type operator. Then the centers of all faces
of the simplex are fixed points of $V$ and
\begin{eqnarray}\label{Thm412}|Fix(V)|\ge
2^m-1,
\end{eqnarray} here as before $|A|$ stands for the cardinality of a set
$A.$
\end{thm}

\begin{pf}
Let $V$ be an $M_1$LV type operator and $x^0=(x_1^0,\dots,x_m^0)$ be
the center of the face $\Gamma_{\alpha},$ i.e.
$$x^0_k=\left\{
\begin{array}{l}
\frac{1}{|\alpha|}, \ \ k\in\alpha \\[2mm]
0, \ \ \ k\notin\alpha.
\end{array}
\right.
$$ where $\alpha\subset I.$

It is clear that $M(x^0)=\alpha.$ According to Theorem
\ref{Volconverge} the trajectory $\{V^nx^0\}_{n=0}^\infty$ converges
to some point $x^{*}$ which is a fixed point of $V.$ Since the face
$\Gamma_\alpha$ is invariant w.r.t. $V$ then $x^{*}\in
\Gamma_\alpha.$ According to Lemma \ref{savemax} we have
\begin{eqnarray}\label{Thm413} M(x^{*})=M(x^0)=\alpha. \end{eqnarray}
which means $x^{*}\in ri\Gamma_\alpha.$ On the other hand, it
follows from \eqref{Thm413} that  all  non null coordinates of
$x^{*}$ are maximal, it means $x^{*}=x^0.$ So, $x^0$ is a fixed
point of $V.$

It is clear that the number of  faces of the simplex is
$$\sum\limits_{i=1}^mC_m^i=2^m-1,$$ so we have \eqref{Thm412}.
\end{pf}

\begin{rem}\label{e1l1}
Note that the operator $V_{\varepsilon,\ell}$ given by
\eqref{VLLEPSILON} was first considered in \cite{UH}, in a
particular case, when $\varepsilon=1,$ $l=1.$ There, it was
established that for every $x^0\in S^{m-1}$ the trajectory
$\{V_{1,1}^nx^0\}_{n=0}^\infty$ starting from any $x^0\in S^{m-1}$
always converges. Since the operator \eqref{VLLEPSILON} is also
$M_1$LV type then according to Theorem \ref{Volconverge} for every
$x^0\in S^{m-1}$ the trajectory
$\{V_{\varepsilon,\ell}^nx^0\}_{n=0}^\infty$ always converges for
all $0<\varepsilon\le1$ and $\ell\in {\mathbb{N}}.$
\end{rem}

Now for the particular operator $V_{\varepsilon,\ell}$ we are going
to find a limit point of the trajectory
$\{V_{\varepsilon,l}^nx^0\}_{n=0}^\infty$.

\begin{obs}\label{omegaVepl}
Let  $x^0\in S^{m-1}$ and $M(x^0)=\{k_1,k_2,\dots,k_r\},$ where
$k_1<k_2<\cdots<k_r,$ then  for the operator $V_{\varepsilon,\ell}$
given by \eqref{VLLEPSILON} one has
$$\lim\limits_{n\rightarrow\infty}V_{\varepsilon,l}^{n}x^0=
\bigg(0,\dots,0,\underbrace{\frac{1}{r}}_{k_1},0,\dots,0,\underbrace{\frac{1}{r}}_{k_2},0,
\dots,0,\underbrace{\frac{1}{r}}_{k_r},0,\dots,0\bigg).$$
\end{obs}

\begin{pf}
Without loss of  generality we may suppose that $x^0\in
\mathrm{ri}S^{m-1},$ otherwise it is enough to consider the
restriction of $V_{\varepsilon,l}$ on the face $\Gamma_{\alpha_0}$
of the simplex, where $\alpha_0=supp(x^0).$ Here $supp(x^0)=\{i\in
I: x_i^0\neq 0\}$.

One can check that all fixed points of the operator
$V_{\varepsilon,\ell}$ are only the centers of faces i.e. if
$\alpha\subset I$ and $\alpha=\{i_1,i_2,\dots,i_r\}$ then
$$Fix(V_{\varepsilon,l})=\bigcup\limits_{\forall\alpha\subset I}
\bigg(0,\dots,0,\underbrace{\frac{1}{|\alpha|}}_{i_1},0,\dots,0,\underbrace{\frac{1}{|\alpha|}}_{i_2},0,\dots,
0,\underbrace{\frac{1}{|\alpha|}}_{i_r},0,\dots,0\bigg)$$

Since the trajectory $\{V_{\varepsilon,\ell}^nx^0\}_{n=0}^\infty$
converges (see Remark \ref{e1l1}) and its the limit point is a fixed
point of $V_{\varepsilon,\ell}$ so we need to show that the limit
point is the center of the face $\Gamma_{M(x^0)}.$

Assume from the contrary that is  the limit point $x^{*}$ of the
trajectory $\{V_{\varepsilon,\ell}^nx^0\}_{n=0}^\infty$ is a center
of another face $\Gamma_\gamma,$ where $\gamma\neq M(x^0).$ Then it
is clear $M(x^{*})=\gamma$ which contradicts to  Lemma
\ref{savemax}.
\end{pf}

It is interesting to know whether is there a quadratic volterrian
operator \eqref{QSO} which is the $M_1$LV type operator.

\begin{thm}\label{VolterianMvol}
A quadratic volterrian operator given by \eqref{QSO} is $M_1$LV type
if and only if one has
$$a_{ki}=0, \ \ \forall \ k,i=\overline{1,m},$$ i.e.
$${\mathcal{VM}}_1\cap {\mathcal{QV}}=\{Id\},$$ here as before  $Id:S^{m-1}\to
S^{m-1}$ is an identity mapping.
\end{thm}

\begin{pf}
\emph{If part.} Let a quadratic volterrian operator given by
\eqref{QSO} be $M_1$LV type. According to Theorem \ref{FixedVol}
each center of all the face of the simplex is a fixed point of
operator given by \eqref{QSO}. In particular, the centers of  all
one dimensional faces
$$x^0=\bigg(0,\dots,0,\underbrace{\frac12}_k,0,\dots,0,\underbrace{\frac12}_i,0,\dots,0\bigg)$$
are fixed points. Then we have the following
$$\left\{
\begin{array}{l}
  (Vx^0)_k=\frac12(1+a_{ki}\frac12)=\frac12 \\[2mm]
  (Vx^0)_i=\frac12(1+a_{ik}\frac12)=\frac12,
\end{array}
\right.
$$
which means $a_{ki}=0$ for all $k,i=\overline{1,m}.$

Only if part is obvious. This completes the proof.
\end{pf}

\begin{prob}
At a moment similar results as  Theorems \ref{FixedVol},
\ref{VolterianMvol} are open for the $M_0$LV type operators.
\end{prob}

\section{$\textbf{f}-$ monotone Lotka-Volterra type operators}

In this section we recall a notion of $\textbf{f}-$monotonicity of
LV type operator and study its properties.

We first  recall \cite{T}, \cite{N} that a mapping
$\textbf{f}:S^{m-1}\to {\mathbb{R}}^m$ is called \emph{monotone} on
the simplex $S^{m-1}$ if for every points $x,y\in S^{m-1}$  the
following condition is satisfied
\begin{eqnarray}\label{51}\langle\textbf{f}(x)-\textbf{f}(y),x-y\rangle\ge0,\end{eqnarray}
where $\langle x,y\rangle$ stands for the standard scalar product in
${\mathbb{R}}^m.$

Observe that for a generating function $\mathbf{f}$ of LV type
operator, its monotonicity can be replaced by
$$
\langle{\mathbf{f}}(x),y\rangle+\langle x,{\mathbf{f}}(y)\rangle\leq
0, \qquad x,y\in S^{m-1}.
$$

\begin{defin}\label{Fmonotone}
An LV type operator \eqref{Volterra} is called $\mathbf{f}-$monotone
if its generating mapping $\textbf{f}:S^{m-1}\rightarrow
{\mathbb{R}}^m$ is monotone on $S^{m-1}$.
\end{defin}

By ${\mathcal{FV}}$ we denote the set  of all $\textbf{f}-$monotone
LV type operators.

\begin{exam}\label{exfmonot}
Now we are going to show that LV type operator
$V_{\varepsilon,\ell}$ (see \eqref{VLLEPSILON}) defined in Example
\ref{M1volterra} is $\textbf{f}-$monotone.
\end{exam}

Indeed, we have \begin{eqnarray*}
    \langle\textbf{f}_{\varepsilon,\ell}(x)-\textbf{f}_{\varepsilon,\ell}(y),x-y\rangle&=&-\langle x,\textbf{f}_{\varepsilon,\ell}(y)\rangle-
    \langle\textbf{f}_{\varepsilon,\ell}(x),y\rangle\\
    &=&-\varepsilon\sum\limits_{k=1}^mx_k\bigg(y_k^\ell-\sum\limits_{i=1}^my_i^{\ell+1}\bigg)
    -\varepsilon\sum\limits_{k=1}^my_k\bigg(x_k^\ell-\sum\limits_{i=1}^mx_i^{\ell+1}\bigg)\\
    &=&-\varepsilon\bigg(\sum\limits_{k=1}^mx_ky_k^\ell+\sum\limits_{k=1}^mx_k^\ell y_k-\sum\limits_{k=1}^mx_k^{\ell+1}-\sum\limits_{k=1}^my_k^{\ell+1}\bigg)\\
     &=&\varepsilon\sum\limits_{k=1}^m(x_k-y_k)(x_k^\ell-y_k^l)\\
    &=&\varepsilon\sum\limits_{k=1}^m(x_k-y_k)^2(x_k^{\ell-1}+x_k^{\ell-2}y_k+...+x_ky_k^{\ell-2}+y_k^{\ell-1})\\
    &\ge&0.
  \end{eqnarray*}
which shows the operator $V_{\varepsilon,\ell}$ is
$\textbf{f}-$monotone.

Note that $\textbf{f}-$monotone LV type operators have the following
properties.

\begin{thm}[\cite{GS}]\label{propf-monot} Let $V$ be a $\emph{\textbf{f}}-$monotone LV
type operator on $S^{m-1}$ then the following assertions hold true:
\begin{itemize}
    \item[(i)] $V$ is a homeomorphism of the simplex.
    \item[(ii)] For any $x^0\in riS^{m-1},$ $Vx^0\neq x^0$ the set  of all limit points $\omega(x^0)$
    of the trajectory $\{V^nx^0\}_{n=0}^\infty$ belongs to the boundary
    $\partial S^{m-1}$ of the simplex. Moreover either
    $|\omega(x^0)|=1$ or $|\omega(x^0)|\ge \aleph_0,$ here as before
    $|\omega(x^0)|$ means the cardinality of the set $\omega(x^0).$
    \item[(iii)] $V$ has no periodic points (except for fixed points) .
\end{itemize}
\end{thm}

\begin{rem}
If $V$ is not $\emph{\textbf{f}}-$monotone then each of the given
statements in Theorem \ref{propf-monot} may not be satisfied. Indeed
let us consider the following examples.
\end{rem}

First we note that in all provided examples below $V$ is not
$\textbf{f}-$monotone.

\begin{exam}
Consider an operator $V$ given by
\begin{equation}\label{53} (Vx)_k=x_k\bigg(x_k^2+3\sum\limits_{i=1}^{k-1}x_i-3\sum\limits_{i,j=1,i<j}^{k-1}x_ix_j\bigg),
\ \ k=\overline{1,m}.
\end{equation}
\end{exam}

Let us first show that the defined operator is LV type.

Since for all $k=\overline{1,m}$ we have
\begin{eqnarray}\label{54}
\sum\limits_{i=1}^{k-1}x_i-\sum\limits_{i,j=1,i<j}^{k-1}x_ix_j=
\sum\limits_{i=1}^{k-2}x_i\left(1-\sum\limits_{j=i+1}^{k-1}x_j\right)+x_{k-1}\ge0,
\end{eqnarray} which implies that $(Vx)_k\ge0$ for any $x\in S^{m-1}$
and $k=\overline{1,m}.$

Denote
$$W_k(x)=\left(\sum\limits_{i=k}^mx_i\right)^3+
3\sum_{i=1}^{k-1}x_i\cdot\left(\sum\limits_{i=k}^mx_i\right)^2+
3\sum\limits_{{i,j=1, \\
i\le j}}^{k-1}x_ix_j\cdot\sum\limits_{i=k}^m x_i, \ \
k=\overline{1,m}.
$$
Then for all $k=\overline{1,m}$ one has

\begin{eqnarray}\label{(*)}
  &&  W_k(x)=(Vx)_k+W_{k+1}(x), \ k=\overline{1,m-1},\\[2mm]
\label{(**)}
 && W_m(x)=(Vx)_m.
\end{eqnarray}

Indeed,
\begin{eqnarray*}
W_k(x)&=&\bigg(x_k+\sum\limits_{i=k+1}^mx_i\bigg)^3+
3\sum\limits_{i=1}^{k-1}x_i\cdot\bigg(x_k+\sum\limits_{i=k+1}^mx_i\bigg)^2\\
& &+3\sum\limits_{{i,j=1, i\le j}}^{k-1}x_ix_j\cdot\bigg(x_k+\sum\limits_{i=k+1}^m x_i\bigg)\\
   &=&x_k\bigg(x_k^2+3\sum\limits_{i=1}^{k-1}x_i\cdot\sum\limits_{i=k}^{m}x_i+3\sum\limits_{{i,j=1,i\le j}}^{k-1}x_ix_j\bigg)+W_{k+1}(x) \\
  &=&x_k\bigg(x_k^2+3\sum\limits_{i=1}^{k-1}x_i-3\bigg(\sum\limits_{i=1}^{k-1}x_i\bigg)^2+3\sum\limits_{{i,j=1,i\le j}}^{k-1}x_ix_j\bigg)+W_{k+1}(x) \\
  &=& (Vx)_k+W_{k+1}(x).
\end{eqnarray*}
Moreover,
\begin{eqnarray*}
  W_m(x) &=&x_m\bigg(x_m^2+3x_m\sum\limits_{i=1}^{m-1}x_i+3\sum\limits_{{i,j=1,i\le j}}^{m-1}x_ix_j\bigg)\\
  &=&x_m\left(x_m^2+3\sum\limits_{i=1}^{m-1}x_i-3\bigg(\sum\limits_{i=1}^{m-1}x_i\bigg)^2+3\sum\limits_{{i,j=1,i\le j}}^{m-1}x_ix_j\right)\\
  &=&x_m\bigg(x_m^2+3\sum\limits_{i=1}^{m-1}x_i-3\sum\limits_{{i,j=1,i<j}}^{m-1}x_ix_j\bigg)\\
  &=&(Vx)_m
\end{eqnarray*}
So from \eqref{(*)} and\eqref{(**)} we find
$$
1=\bigg(\sum\limits_{i=1}^mx_i\bigg)^3=W_1(x)=(Vx)_1+W_2(x)=\cdots=\sum\limits_{i=1}^{m}(Vx)_i.
$$
Therefore, the operator given by \eqref{53} is LV type.

Now we show that the operator $V$ is injective. Indeed, let $x,y\in
S^{m-1}$ and $x\neq y.$ Then there exists $k_0\in I$ such that
$$x_{k_0}\neq y_{k_0}, \ \ \  x_i=y_i , \ \ \forall i=\overline{1,k_0-1}.$$
From \eqref{53} we find that
$$(Vx)_i=(Vy)_i, \ \ i=\overline{1,k_0-1}.$$
The equality \eqref{54} with $x_i=y_i$ for all
$i=\overline{1,k_0-1}$ yields that
$$C:=\sum\limits_{i=1}^{k_0-1}x_i-\sum\limits_{{i,j=1,i<j}}^{k_0-1}x_ix_j=
\sum\limits_{i=1}^{k_0-1}y_i-\sum\limits_{{i,j=1,i<j}}^{k_0-1}y_iy_j\ge
0.$$ Now consider a function
$$g(t)=t^3+3C\cdot t, \ \ t\in [0,1],$$
which is strictly increasing on the segment $[0,1].$ Therefore, for
$x_{k_0}\neq y_{k_0}$ we get
$$(Vx)_{k_0}=g(x_{k_0})\neq g(y_{k_0})=(Vy)_{k_0},$$
which means  $Vx\neq Vy$.

According to Corollary \ref{gomeom} we can conclude that the
operator $V$ is a homeomorphism.

Let us show that \eqref{51} is not satisfied. From \eqref{53} we
find the corresponding generating function $\mathbf{f}$ of $V$ as
follows
$$({\mathbf{f}}(x))_k=x_k^2+3\sum\limits_{i=1}^{k-1}x_i-3\sum\limits_{{i,j=1,i<j}}^{k-1}x_ix_j-1, \ \ k=\overline{1,m}.$$
For $x^0=(1,0,\dots,0)$ and $y^0=(0,1,0,\dots,0)$ one has
$$
    \langle\textbf{f}(x^0)-\textbf{f}(y^0),x^0-y^0\rangle=-1<0,
$$ which implies that $V$ is not  $\textbf{f}-$monotone.

\begin{rem}
Thus, the provided example shows that a LV type operator
\eqref{Volterra} to be a homeomorphism  the
$\textbf{f}-$monotonicity is a sufficient condition.
\end{rem}

Let us provide another example for a LV type operator
\eqref{Volterra} which it is not $\textbf{f}-$monotone and does not
satisfy the assertion $(ii)$ of Theorem \ref{propf-monot}.

\begin{exam} Let us consider $M_0$LV type operator $W_{\varepsilon,\ell}$ defined
in Example \ref{M0Volterra}.
\end{exam}

Using the same argument of Example \ref{exfmonot} we can show that
$W_{\varepsilon,\ell}$ is  not $\textbf{f}-$monotone, i.e. one has
if $x,y\in S^{m-1}$ with $x\neq y$ then
$$\langle\textbf{f}_\varepsilon(x)-\textbf{f}_\varepsilon(y),x-y\rangle<0.$$

According to Theorem \ref{Volconverge} for any $x^0\in S^{m-1}$ the
trajectory $\{W_{\varepsilon,\ell}^nx^0\}_{n=0}^\infty$  converges.
Now we find its limit point.

\begin{obs}\label{omegaWEPL}
Let $x^0\in ri\Gamma_\alpha,$ where $\alpha$ is any subset of $I,$
then a limit point $\omega(x^0)$ of the trajectory
$\{W_{\varepsilon,\ell}^nx^0\}_{n=0}^\infty$  is the center of the
face $\Gamma_\alpha.$
\end{obs}

\begin{pf}
It is obvious that the fixed points of operator
$W_{\varepsilon,\ell}$  (see \eqref{WEPSILONL} ) are only the
centers of all face of the simplex i.e. if $\alpha\subset I$ and
$\alpha=\{i_1,i_2,\cdots,i_r\}$ then
$$Fix(W_{\varepsilon,l})=\bigcup\limits_{\forall\alpha\subset I}
\bigg(0,\dots,0,\underbrace{\frac{1}{|\alpha|}}_{i_1},0,\dots,0,
\underbrace{\frac{1}{|\alpha|}}_{i_2},0,\dots,0,
\underbrace{\frac{1}{|\alpha|}}_{i_r},0,\dots,0\bigg)$$

Since the trajectory $\{W_{\varepsilon,\ell}^nx^0\}_{n=0}^\infty$
converges  and the limit point of the trajectory is a fixed point of
$W_{\varepsilon,\ell}$  so we need to show that the limit point is
the center of the  face $\Gamma_\alpha.$

Without loss of  generality we may suppose that $x^0\in
\mathrm{ri}S^{m-1}.$

From \eqref{WEPSILONL} one gets
\begin{eqnarray}\label{(***)}
(W_{\varepsilon,l}x^0)_k=x_k^0\left(1+\varepsilon\sum\limits_{i=1}^m\left((x_i^0)^l-(x_k^0)^l\right)x_i^0\right),
\ \ k=\overline{1,m}. \end{eqnarray}

Let $$m(x^0)=\{k\in I: x_k^0=\min\limits_{i=\overline{1,m}}x_i^0\}$$
and take $k\in m(x^0).$ Then it follows from \eqref{(***)} that
\begin{eqnarray}\label{56}0<x_k^0\le (W_{\varepsilon,l} x^0)_k\le
(W_{\varepsilon,\ell}^{2} x^0)_k\le \cdots.   \end{eqnarray} which
means  that
$$\lim\limits_{n\rightarrow\infty}(W_{\varepsilon,\ell}^{n} x^0)_k>0, \ \ k\in m(x^0),$$
i.e. the minimal coordinate of the limit point is positive. Then the
limit point belongs to  $riS^{m-1}$.  It is clear that the interior
fixed point of $W_{\varepsilon,\ell}$ is only the center of the
simplex, so we obtain
$$\lim\limits_{n\rightarrow\infty}W_{\varepsilon,\ell}^{n}x^0=\left(\frac1m,\frac1m,\cdots,\frac1m\right).$$
\end{pf}

From Observation \ref{omegaWEPL} we immediately find that $(ii)$ of
Theorem \ref{propf-monot} is not satisfied for
$W_{\varepsilon,\ell}.$

\begin{rem} The provided example shows that
${\mathcal{VM}}_0\neq {\mathcal{FV}},$ since
$W_{\varepsilon,l}\notin {\mathcal{FV}}.$
\end{rem}

Now we are going to show another example of LV type operator which
is not $\textbf{f}-$monotone and does not satisfy assertion $(iii)$
of Theorem \ref{propf-monot}.

\begin{exam}\label{notinjecperiodic}
Consider one dimensional simplex $S^1$ with the following
decomposition
$$S^1=T_1\bigcup T_2\bigcup T_3 \bigcup T_4\bigcup T_5,$$
where $$T_1=\left\{x\in S^1: 0\le x_1\le\frac{9}{30}\right\},
T_2=\left\{x\in S^1: \frac{9}{30}\le x_1\le\frac{11}{30}\right\},$$
$$T_3=\left\{x\in S^1: \frac{11}{30}\le x_1\le\frac{19}{30}\right\},$$
$$T_4=\left\{x\in S^1: \frac{19}{30}\le x_1\le\frac{21}{30}\right\}, T_5=\left\{x\in S^1: \frac{21}{30}\le x_1\le 1\right\}.$$
Define an LV type operator $V:S^1\rightarrow S^1$  by
\begin{eqnarray}\label{57}
\left\{%
\begin{array}{ll}
    (Vx)_1=x_1(1+f_1(x_1,x_2)) \\
    (Vx)_2=x_2(1+f_2(x_1,x_2)) \\
\end{array}%
\right. \end{eqnarray} here
$$f_1(x_1,x_2)=\begin{cases}
    2 & \text{if \ \  $x\in T_1$,} \\
    \dfrac{9}{10x_1}-1 & \text{if \ \ $x\in T_2$,} \\
    \dfrac{4x_2-2}{x_1} & \text{if \ \  $x\in T_3$,} \\
    \dfrac{1}{10x_1}-1  & \text{if \ \ $x\in T_4$,} \\
    -\dfrac{2x_2}{x_1} & \text{if \ \  $x\in T_5$,}
\end{cases}
\ \ \ \ f_2(x_1,x_2)=\begin{cases}
    -\dfrac{2x_1}{x_2} & \text{if \ \  $x\in T_1$,} \\
    \dfrac{1}{10x_2}-1 & \text{if \ \ $x\in T_2$,} \\
    \dfrac{4x_1-2}{x_2} & \text{if \ \  $x\in T_3$,} \\
    \dfrac{9}{10x_2}-1  & \text{if \ \ $x\in T_4$,} \\
    2 & \text{if \ \  $x\in T_5$.}
\end{cases}
$$
\end{exam}
One can see that
$$V(T_2)=\left(\frac{9}{10},\frac{1}{10}\right), \ \  V(T_4)=\left(\frac{1}{10},\frac{9}{10}\right),$$
i.e. the operator given by \eqref{57} is not injective (moreover it
is not a homeomorphism, see Corollary \ref{gomeom}) this implies
that $V$ is not $\textbf{f}-$monotone (see (i) of Theorem
\ref{propf-monot}). But $V$ has 2-periodic points in $S^1$. Indeed,
$$V^2\left(\dfrac14,\dfrac34\right)=\left(\dfrac14,\dfrac34\right), \ \
V^2\left(\dfrac34,\dfrac14\right)=\left(\dfrac34,\dfrac14\right).$$

This establishes that the condition $(iii)$ of Theorem
\ref{propf-monot} is not satisfied.\\

\begin{rem} Finally we have to strass that the set ${\mathcal{VM}}_1$ of
$M_1$LV type operators and the set ${\mathcal{FV}}$ of
$\textbf{f}-$monotone LV type operators does not coincide i.e.
${\mathcal{VM}}_1\neq {\mathcal{FV}}.$ Indeed, consider the
following
\end{rem}

\begin{exam}
Consider again one dimensional simplex $S^1$ with its following
decomposition
$$S^1=T_1\bigcup T_2\bigcup T_3 \bigcup T_4,$$
where $$T_1=\left\{x\in S^1: 0\le x_1\le\frac{1}{3}\right\},
T_2=\left\{x\in S^1: \frac{1}{3}\le x_1\le\frac{5}{12}\right\},$$
$$T_3=\left\{x\in S^1: \frac{5}{12}\le x_1\le\frac{1}{2}\right\}, T_4=\left\{x\in S^1: \frac{1}{2}\le x_1\le 1\right\}.$$
Define an LV type operator $V:S^1\rightarrow S^1$  by
\begin{eqnarray}\label{58}
\left\{%
\begin{array}{ll}
    (Vx)_1=x_1(1+f_1(x_1,x_2)) \\
    (Vx)_2=x_2(1+f_2(x_1,x_2)) \\
\end{array}%
\right. \end{eqnarray} here
$$
f_1(x_1,x_2)=\begin{cases}
    0 & \text{if \ \  $x\in T_1$,} \\
    \dfrac{x_2-2x_1}{3x_1} & \text{if \ \ $x\in T_2$,} \\
    \dfrac{x_1-x_2}{2x_1} & \text{if \ \  $x\in T_3$,} \\
    0  & \text{if \ \ $x\in T_4$,}
\end{cases}
\ \ \ \ f_2(x_1,x_2)=\begin{cases}
    0 & \text{if \ \  $x\in T_1$,} \\
    \dfrac{2x_1-x_2}{3x_2} & \text{if \ \ $x\in T_2$,} \\
    \dfrac{x_2-x_1}{2x_2} & \text{if \ \  $x\in T_3$,} \\
    0  & \text{if \ \ $x\in T_4$,}
\end{cases}$$
\end{exam}

One can show that the operator given by \eqref{58} is $M_1$LV type.
It is easy to see that $V(T_2)=(\frac{1}{3},\frac{2}{3}),$ i.e. the
operator is not injective (moreover it is not a homeomorphism, see
Corollary \ref{gomeom}) which means it is not $\textbf{f}-$monotone
(see (i) of Theorem \ref{propf-monot}).

\section*{Acknowledgement}

The authors acknowledge Research Endowment Grant B (EDW B 0905-303)
of IIUM and the MOSTI grant 01-01-08-SF0079.

\end{document}